\documentclass[11pt]{article}

\usepackage{amsmath,amssymb}

\newtheorem{theorem}{Theorem}
\newtheorem{ass}{Theorem}[section]
\newtheorem{prop}[ass]{Proposition}
\newtheorem{lemma}[ass]{Lemma}

\newtheorem{corollary}[ass]{Corollary}

\newtheorem{conj}{Conjecture}
\newtheorem{prob}{Problem}
\newtheorem{definition}[ass]{Definition}

\newcommand{\qed}{\hspace*{\fill} \rule{7pt}{7pt}}
\newcommand{\Proof}{\noindent{\bf Proof}\ \ }

\begin{document}

\title{A tournament approach to pattern avoiding matrices}

\author{
Asaf Shapira
\thanks{School of Mathematics, Tel-Aviv University, Tel-Aviv, Israel 69978. Email: {\tt asafico@tau.ac.il}. Supported in part by ISF Grant 224/11 and a Marie-Curie CIG Grant 303320.}
\and
Raphael Yuster
\thanks{Department of Mathematics, University of Haifa, Haifa
31905, Israel. Email: raphy@math.haifa.ac.il}
}

\date{}

\maketitle

\setcounter{page}{1}

\begin{abstract}

We consider the following Tur\'an-type problem: given a fixed tournament $H$,
what is the least integer $t=t(n,H)$ so that adding $t$ edges to any $n$-vertex tournament,
results in a digraph containing a copy of $H$. Similarly, what is the least integer $t=t(T_n,H)$
so that adding $t$ edges to the $n$-vertex transitive tournament,
results in a digraph containing a copy of $H$. Besides proving several results on these problems,
our main contributions are the following:

\begin{itemize}

\item Pach and Tardos conjectured that if $M$ is an acyclic $0/1$ matrix, then any
$n \times n$ matrix with $n(\log n)^{O(1)}$ entries equal to $1$ contains the pattern $M$. We show that this conjecture
is equivalent to the assertion that $t(T_n,H)=n(\log n)^{O(1)}$ if and only if $H$ belongs to a certain (natural) family of tournaments.

\item We propose an approach for determining if $t(n,H)=n(\log n)^{O(1)}$.
This approach combines expansion in sparse graphs, together with certain structural characterizations
of $H$-free tournaments.
\end{itemize}

Our result opens the door for using structural graph theoretic tools in order to settle the Pach-Tardos conjecture.


\end{abstract}

\section{Introduction}

Tur\'an-type problems are among the most basic and most well studied problems in extremal
combinatorics. In the setting of graphs, Tur\'an's problem asks, given a fixed graph $H$ and an integer $n$,
what is the least integer $m=\mbox{ex}(n,H)$, so that every $n$-vertex graph with at least $m$ edges contains a copy of $H$.
Tur\'an's classical theorem \cite{turan-1941} solved the problem when $H$ is a complete graph. Since then, Tur\'an-type problems
have been extensively studied in numerous settings. In this paper, we propose to study a natural Tur\'an-type problem in
the setting of tournaments. While part of our motivation was a recent surge in the study of tournaments, our main motivation
was a surprising connection between this new problem and a classical Tur\'an-type problem is the setting of matrices.

\subsection{A Tur\'an-type problem for tournaments}

A {\em tournament} $G=(V,E)$ is a digraph such that for every pair of distinct vertices $u,v$, the edge set $E$ contains exactly one edge with
ends $\{u,v\}$. In other words, either $(u,v) \in E$ or $(v,u) \in E$ is present, but not both.
If $G=(V,E)$ is a tournament, we say that $X \subseteq V$ is {\em transitive} if the sub-tournament $G[X]$ induced on $X$ has
no directed cycle. We denote by $T_n$ the (unique) transitive tournament on $n$ vertices.
A pair $\{u,v\}$ in a digraph is {\em bidirectional} if both $(u,v)$ and $(v,u)$ are edges.
Notice that if we add $t$ new edges to a tournament, the resulting digraph has $t$ bidirectional pairs, or, equivalently,
$t$ cycles of length $2$. Tournaments having some bidirectional edges are known in the literature as {\em semi-complete digraphs}.

Let $H$ be any given tournament. We say that a digraph is {\em $H$-free} if it contains no copy of $H$ as a subgraph.
Clearly, if $G$ is the complete digraph\footnote{The complete digraph on $n$ vertices is the digraph containing all possible $2\binom{n}{2}$ edges.} on $n \ge |V(H)|$ vertices, then $G$ is not $H$-free.
On the other hand, if $H$ is any non-transitive tournament, then any transitive tournament is $H$-free.
If we take any tournament $G$ on $n \ge |V(H)|$ vertices and start adding edges to $G$ (thereby obtaining a semi-complete digraph)
we will, at some point, obtain a digraph which is not $H$-free. This motivates the following problems.

\begin{prob}\label{pr:TTP}
Given a tournament $H$ and an integer $n$, determine the smallest integer $t=t(n,H)$ so that
adding $t$ edges to any $n$-vertex tournament $G$, results in a digraph containing a copy of $H$.
\end{prob}
\begin{prob}\label{pr:TTP-TN}
Given a tournament $H$ and an integer $n$, determine the smallest integer $t=t(T_n,H)$ so that
adding $t$ edges to $T_n$ results in a digraph containing a copy of $H$.
\end{prob}
Notice that $t(n,H)$ is the minimum integer such that any semi-complete digraph with at least $\binom n2 +t(n,H)$ edges has a copy of $H$.
We always assume that $n \ge |V(H)|$ (as otherwise $t(n,H)$ and $t(T_n,H)$ are undefined).
Observe also that $t(T_n,H) \leq t(n,H)$ for every $H$
and that for any non-transitive tournament $H$ we  have $t(n,H) \ge n/2$ as one can make bidirectional the edges of a matching in $T_n$ without introducing a cycle of length $3$. We will be mainly interested in the asymptotic behaviors of $t(n,H)$ and $t(T_n,H)$ for a fixed $H$ and $n \rightarrow \infty$.

\subsection{Some preliminary observations and a conjecture}\label{subsec:prelim}

As it turns, similar to the case of Tur\'an's problem for graphs, there are certain
conditions under which Problems \ref{pr:TTP} and \ref{pr:TTP-TN} are quite easy to answer.
We devote this subsection to these cases, deferring the more interesting cases to the following ones.
We start by recalling the notion of a coloring of a tournament, which was introduced in \cite{BCCFLSST-2013}.
For a positive integer $k$, a {\em $k$-coloring} of a tournament $G=(V,E)$ is a partition of $V$ into $k$
parts, where each part induces a transitive sub-tournament. The {\em chromatic number} $\chi(G)$
of a tournament $G$ is the minimum $k$ such that $G$ admits a $k$-coloring.
By definition, $\chi(G)=1$ if and only if $G$ is transitive.

It is not difficult to determine $t(n,H)$ and $t(T_n,H)$ asymptotically whenever $\chi(H) \ge 3$.
It turns out that when $\chi(H)=r \ge 3$, both $t(n,H)$ and $t(T_n,H)$ are very close
to the Tur\'an number of (undirected) graphs with chromatic number $r$.
Recall that for an undirected graph $U$, we take $\mbox{ex}(n,U)$ to be the smallest integer so that every $n$-vertex graph
with at least $\mbox{ex}(n,U)$ edges has a copy of $H$.
The well-known Erd\H{o}s-Stone-Simonovits Theorem \cite{ES-1946} asserts that if $\chi(U)=r$, then $\mbox{ex}(n,U)=(1-\frac{1}{r-1}+o_n(1))\binom{n}{2}$.
A similar value holds for $t(n,H)$ and $t(T_n,H)$.

\begin{prop}\label{p:non-2-chrom}
Let $H$ be a tournament with $\chi(H)=r$. Then
$$
\left(1-\frac{1}{r-1}\right)\binom{n}{2} \le t(T_n,H) \le t(n,H)=\left(1-\frac{1}{r-1}+o_n(1)\right)\binom{n}{2}\;.
$$
\end{prop}

Proposition \ref{p:non-2-chrom} leaves open the asymptotic values of $t(n,H)$ and $t(T_n,H)$ whenever $\chi(H)=2$.
Therefore, all remaining results in this paper have to do with the case $\chi(H)=2$.

Our next observation is that, similar to the case of undirected graphs, whenever $\chi(H)=2$, we have $t(n,H) < n^{2-c}$ for some $c>0$.
To state this result we will need the following definition: for a $2$-chromatic tournament $H$, let $s(H)$ denote the smallest integer such that there is a $2$-coloring of $H$ with a color class of size $s$.

\begin{prop}\label{p:2-chrom-upper}
If $\chi(H) \leq 2$, then $t(n,H)=O(n^{2-1/2^{s(H)-1}})$. In particular, if $s(H)=1$, then $t(n,H)=O(n)$.
\end{prop}

As we will see later, one can obtain significantly stronger bounds for certain tournaments $H$.

Now that we know that either $t(n,H)=\Theta(n^2)$ or  $t(n,H) = O(n^{2-c})$, it is natural
to ask when do we have $t(n,H) = \Omega(n^{1+c})$ for some $c=c(H)>0$. To this end, we need to introduce an important definition.

\begin{definition}[Tournament Forest]\label{TourForest}
A tournament $H$ is said to be a {\em tournament forest} if it is possible to partition its vertex set into
two transitive sets\footnote{Note that this means that if $H$ is a tournament forest, then $\chi(H) \leq 2$.} $L,R$ so that if we define an undirected bipartite graph $U$ on the vertex sets $L$ and $R$ by connecting $\ell \in L$ to $r \in R$ if and only if $H$ has an edge pointing from
$r$ to $\ell$, then $U$ is a forest.
\end{definition}

The next proposition supplies a sufficient condition for $H$ to satisfy $t(n,H) \geq \Omega(n^{1+c})$.

\begin{prop}\label{p:non-strong}
If $H$ is not a forest, then there is $c=c(H)>0$ so that $$t(n,H) \geq t(T_n,H) \geq \Omega(n^{1+c})\;.$$
Moreover, for every $\epsilon > 0$, there is a tournament $H$ with $\chi(H)=2$ satisfying $t(T_n,H)=\Omega(n^{2-\epsilon})$.
\end{prop}

Recapping, Propositions \ref{p:non-2-chrom}, \ref{p:2-chrom-upper} and \ref{p:non-strong} imply
that if $\chi(H)>2$, then $t(n,H)=\Theta(n^2)$ and if $\chi(H)=2$ and $H$ is not a forest, then $n^{1+c'}\leq t(n,H) \leq n^{2-c}$.
In other words, up to now the situation is similar to the case of Tur\'an's problem for undirected graphs.
The only tournaments not covered by the above results are thus the tournament forests. We conjecture that the following
result holds.

\begin{conj}\label{c:MainConj}
If $H$ is a tournament forest, then $t(T_n,H)=n(\log n)^{O(1)}$.
\end{conj}

Note that for an undirected graph $U$ satisfying $\chi(U) \leq 2$, then as for the growth rate of the Tur\'an function we either have $\mbox{ex}(n,U)=O(n)$ or $\mbox{ex}(n,U) = \Omega(n^{1+c})$. More precisely, $\mbox{ex}(n,U)=O(n)$ whenever $U$ is a forest, and $\mbox{ex}(n,U) = \Omega(n^{1+c})$ otherwise.
We will see later that there are tournament forests satisfying $t(T_n,H) = \omega(n \log n)$, and hence
assuming Conjecture \ref{c:MainConj} we get that there are tournaments satisfying $t(T_n,H)=n(\log n)^{\Theta(1)}$.
We will also see that as opposed to the case of graphs, it is probably quite hard to tell if all tournament
forests satisfy $t(T_n,H) = n(\log n)^{O(1)}$.

\subsection{The Pach-Tardos Conjecture}\label{subsec:PTC}

We now turn to discuss the first main result of this paper, which establishes a surprising relation
between Conjecture \ref{c:MainConj} and a conjecture of Pach and Tardos \cite{PT-2006}.
To this end, we need to recall another famous Tur\'an-type problem, this time involving $0/1$ matrices.
Given two $0/1$ matrices $M$ and $A$, where $M$ is $k \times k$ and $A$ is $n \times n$, we say that $A$ contains the pattern
$M$, if $A$ has $k$ rows $r_1 < r_2 < \cdots < r_k$ and $k$ columns $c_1 < c_2 \cdots < c_k$ so that $A_{c_i,r_j}=1$ whenever $M_{i,j}=1$.
The corresponding Tur\'an problem then asks what is the least integer $m=\mbox{ex}(n,M)$ so that every $n \times n$ matrix $A$ with $m$
entries equal to $1$ contains the pattern $M$. This problem was first introduced by F{\"u}redi and Hajnal \cite{FH-1992},
who showed that certain extremal problems can be reduced to the study of $\mbox{ex}(n,M)$.
This problem received a lot of attention in the past two decades, due to its relation to many other problems. Perhaps the most famous result in this
area is due to Marcus and Tardos \cite{MT-2004} who proved that if $M$ is a permutation matrix, then $\mbox{ex}(n,M)=O(n)$, thus
establishing the famous Stanley-Wilf conjecture. The study of $\mbox{ex}(n,M)$ has also found some surprising applications in theoretical
computer science \cite{pettie-2010}.

If one thinks of $M$ and $A$ as being the adjacency matrices
of two {\em ordered} bipartite graphs, then $A$ contains $M$ if and only if the ordered graph defined by $A$ contains an ordered copy of the
ordered graph defined by $M$. Note that with this interpretation in mind, if $M$ is a permutation matrix, then it defines
an ordered bipartite graph which forms a matching. Let us then say that that a $0/1$ matrix is a forest, if the bipartite
graph it defines is a forest. In an attempt to unify several results concerning the matrix Tur\'an problem mentioned above,
Pach and Tardos raised the following conjecture.

\begin{conj}[Pach-Tardos \cite{PT-2006}]\label{c:PTC}
If the matrix $M$ is a forest then $\mbox{ex}(n,M)=n(\log n)^{O(1)}$.
\end{conj}

We obtain the following.

\begin{theorem}\label{t:ConjEquiv}
Conjectures \ref{c:MainConj} and \ref{c:PTC} are equivalent.
\end{theorem}

Note that Conjecture \ref{c:PTC} is concerned with $0/1$ matrices or equivalently {\em ordered/labeled} bipartite graphs.
These are inherently different from ``usual'' graphs
in which we do not think of the vertices as having labels. The importance of Theorem \ref{t:ConjEquiv} thus lies in showing
that Conjecture \ref{c:PTC} is actually equivalent to a problem involving unlabeled/unordered digraphs. Hence, one can now try
to apply standard graph theoretic tools in order to resolve it. We elaborate on this in the next subsection.

\subsection{A tournament approach for the Pach-Tardos Conjecture}\label{subsec:PTCapp}

We now describe an approach for resolving Conjecture \ref{c:PTC} which relies heavily on both
structural and extremal graph theoretic tools. Let us say that a (usual undirected) graph $U=(V,E)$ is a
$d$-expander if for every $S \subseteq V$ of size  $1 \leq |S| \leq |V|/2$ there are
at least $d|S|$ edges connecting $S$ to $V \setminus S$. Let us say that a semi-complete digraph $G=(V,E)$ is a $d$-expander
if the undirected graph $U=(V,E')$ in which $(i,j) \in E'$ if and only if $\{i,j\}$ is a bidirectional edge in $G$, is a $d$-expander.

Given a semi-complete digraph $F$ let ${\cal D}_F$ be the distribution obtained by independently and uniformly
picking, for each bidirectional edge $\{i,j\}$ of $F$, either the edge $(i,j)$ or the edge $(j,i)$.
Note that ${\cal D}_F$ is a distribution over tournaments.

\begin{definition}[Sparse Family]\label{def:sparse} A family of tournaments ${\cal Q}$ is {\em sparse} if
there is a constant $b$ so that the following holds for every $n$-vertex semi-complete digraph $F$ which is a
$b\log^b n$-expander:
\begin{equation}\label{eq:sparse}
\mathbb{P}_{G \sim {\cal D}_F}[G \in {\cal Q} ] < 1/2\;.
\end{equation}
\end{definition}

Given a tournament $G=(V,E)$, a {\em homogeneous set} of $G$ is a subset $X \subseteq V$ such that for
all $v \in V \setminus X$, either $(v,x) \in E$ for all $x \in X$ or $(x,v) \in E$ for all $x \in X$.
A homogeneous set $X$ is nontrivial if $1 < |X| < |V|$, otherwise it is trivial. A tournament is {\em prime} if all of its homogeneous sets are trivial.
Prime tournaments form an important subclass of tournaments. For example, Alon, Pach and Solymosi \cite{APS-2001} have shown that
the famous Erd\H{o}s-Hajnal conjecture is equivalent to the analogous conjecture for prime tournaments. For other related recent results
see \cite{liu-2014,Choromanski-2014} and their references.

\begin{definition}[$H$-useful]\label{def:useful}
A family of tournaments ${\cal Q}_H$ is $H$-useful if every prime tournament which is $H$-free belongs to ${\cal Q}_H$.
\end{definition}

\begin{theorem}\label{t:sparse} If there exists a sparse $H$-useful family,
then $t(n,H)=n(\log n)^{O(1)}$.
\end{theorem}

We thus see that if every tournament forest $H$ has a sparse $H$-useful family of tournaments,
then Conjecture \ref{c:MainConj} holds. Observe that the approach of Theorem \ref{t:sparse} would
actually yield a near linear bound for $t(n,H)$ and not only for $t(T_n,H)$.
In particular, Theorems \ref{t:ConjEquiv} and \ref{t:sparse} have the following corollary,
summarizing our approach to the Pach-Tardos Conjecture.

\begin{corollary}\label{c:NewApp}
If every tournament forest $H$ has a sparse $H$-useful family, then Conjecture \ref{c:PTC} holds.
\end{corollary}

Informally speaking, what the above corollary shows is that in order to resolve Conjecture \ref{c:PTC} it is
enough to understand the structure of prime $H$-free tournaments when $H$ is a forest.
We note that understanding the structure of $H$-free tournaments for arbitrary $H$ is most like
hopeless, so it is important that one needs to only consider the case of $H$ being a forest and
the $H$-free graph being prime.

It is natural to ask at this point if there are non-trivial tournaments $H$ for which one can provably show that there exists
a sparse $H$-useful family. Indeed, there are several such tournaments. See for example \cite{latka-2003,liu-2014} which prove the existence of nontrivial $H$-useful families
for certain tournaments (all of which can be shown to be sparse).
We focus on one of them in particular, as for this particular tournament, the upper bound we obtain using Theorem \ref{t:sparse} is better
than any other bound we can achieve using other results in this paper.

Consider the tournament $U_5$
defined as follows. Let $T_5$ be the transitive tournament on $\{1,2,3,4,5\}$ in which $(i,j)$ is an edge whenever $i < j$.
Then $U_5$ is the tournament obtained from $T_5$ by reversing the direction of the edges $(2,5)$ and $(1,4)$.
For an odd integer $n$ we define $C_n$ to be the tournament on vertices $\{1,\ldots,n\}$ in which $(i,j)$ is an edge
if and only if $j-i ~(\mbox{mod } n) \in \{1,\ldots,(n-1)/2\}$.
We note that $U_5$ has several interesting properties, as shown in \cite{BCCFLSST-2013,liu-2014}.
Liu has recently obtained the following structural result.

\begin{theorem}[Liu \cite{liu-2014}]\label{t:prime}
If $G=(V,E)$ is a prime tournament that is $U_5$-free, then $G$ satisfies one of the following:
\begin{enumerate}
\item $G$ is isomorphic to $C_n$ for some odd $n$.
\item $V$ can be partitioned into $3$ sets $V_1,V_2,V_3$ so that for every $1 \leq i < j \leq 3$ the tournament $G[V_i \cup V_j]$ is transitive.
\end{enumerate}
\end{theorem}

Let ${\cal Q}_5$ be the family of tournaments that satisfy either condition $1$ or $2$ of Theorem \ref{t:prime}.
Then ${\cal Q}_5$ is $U_5$-useful. We further prove the following.

\begin{theorem}\label{t:u5useful}
${\cal Q}_5$ is a sparse family. Hence, $t(n,U_5)=n(\log n)^{O(1)}$.
\end{theorem}

\subsection{Further results regarding $t(n,H)$ and $t(T_n,H)$}\label{subsec:add}

As a step towards Conjecture \ref{c:MainConj}, it is natural to look for sufficient conditions
that would guarantee that $t(T_n,H)$ is indeed of order $n(\log n)^{O(1)}$.
In this subsection we describe several cases in which we can verify this conjecture, but we would first like to raise
the following interesting problem.

\begin{prob}
Is it true that for any $H$ we have $t(n,H) \sim t(T_n,H)$\;?
\end{prob}

Recall that we trivially have $t(T_n,H) \leq t(n,H)$ thus the problem lies in proving the other direction.
As of now, we do not have any example ruling out even the possibility that $t(n,H) = t(T_n,H)$ for all $n$ sufficiently large.

For a given linear order of the vertices of tournament $H=(V,E)$, we say that $(i,j)$ is a {\em back edge} if $j$ appears before $i$ in the order.
The set of back edges of a given ordering of $H$ forms a {\em back-edge graph} which one can view as an undirected graph.
The {\em minimum feedback edge set} number of a tournament $H$, denoted by $\beta(H)$ is the smallest number of edges in a back-edge graph of $H$.
Observe that if $\beta(H) \ge |V(H)|$, then any back-edge graph has a cycle.
We say that a tournament is a {\em weak forest} if it has an acyclic back-edge graph. Notice that weak forest tournaments are trivially $2$-chromatic.
Observe that every tournament forest (as defined in Definition \ref{TourForest}) is also a weak forest since we can take
the order which puts all vertices of $L$ before all vertices of $R$.
We note however, that being a forest is a stronger requirement.
For example, consider the tournament $\Delta_k$ having $k$ strongly connected components, each of which is a $3$-cycle $C_3$.
It is not hard to see that this tournament is a weak forest but not a forest for all $k \ge 3$.

So for which tournaments can we (unconditionally) show that $t(T_n,H)=n(\log n)^{O(1)}$?
Note that Proposition \ref{p:2-chrom-upper} shows that this is the case when $s(H)=1$.
Observe that $s(H)=1$ if and only if $H$ is a weak forest that has a back-edge graph consisting of a star and some isolated
vertices. For this reason we call such tournaments {\em star tournaments}.
Observe that star tournaments are actually tournament forests as one may take the root of the star as the side $L$ and the remaining vertices
as the side $R$. We now turn to show that there are other cases in which we can prove linear bounds.

The following result shows that if a tournament has a back-edge graph that is very small,
then indeed $t(T_n,H)$ is linear.
\begin{theorem}\label{t:beta-2}
Any tournament $H$ with $\beta(H) \le 2$ satisfies $t(T_n,H)=O(n)$.
\end{theorem}

We note that Theorem \ref{t:beta-2} cannot be extended to hold whenever $\beta(H) \le 3$.
Indeed the tournament $\Delta_3$ mentioned earlier satisfies $\beta(\Delta_3)=3$  and one can
check that it is not a forest, so Proposition \ref{p:non-strong} implies that $t(T_n,\Delta_3)> n^{1+c}$.

Let us end with two more observation. Recall that Proposition \ref{p:2-chrom-upper} implies that
$t(n,H)=n(\log n)^{O(1)}$ whenever $H$ is a star tournament. Theorem \ref{t:beta-2} shows that $t(T_n,H)=n(\log n)^{O(1)}$
also for certain tournaments forests that are not stars. To see that there are also non-star
tournaments satisfying the stronger condition $t(n,H)=n(\log n)^{O(1)}$ observe that
that the tournament $U_5$ defined in the previous subsection is not a star tournament.
Hence, Theorem \ref{t:u5useful} shows that there are non-star tournaments for which $t(n,H)$ is close to linear.
Second, as we mentioned at the end of Subsection \ref{subsec:prelim},
the study of $t(n,H)$ and $t(T_n,H)$ when $H$ is a forest appears much harder than the corresponding problem
for graphs. In particular, while for undirected graphs we have $\mbox{ex}(n,H)=O(n)$ whenever $H$ is a forest,
the following result shows that this is not the case for tournaments.

\begin{theorem}\label{t:non-linear}
There are forest tournaments satisfying $t(T_n,H) = \omega(n \log n)$.
\end{theorem}

\subsection{Paper organization}

In Section \ref{sec:prelim} we prove the preliminary observations given in Propositions \ref{p:non-2-chrom}, \ref{p:2-chrom-upper} and \ref{p:non-strong}.
While the proofs are quite routine, we also prove in the same section more refined versions of these results which
apply more ideas. The first main result of this paper, stated in Theorems \ref{t:ConjEquiv}, is proved in Section \ref{sec:equiv}.
While one direction of this equivalence is simple, the other direction requires a somewhat delicate reduction.
In the same section we prove Theorem \ref{t:non-linear}.
The second main result of this paper, stated in Theorem \ref{t:sparse}, is proved
in Section \ref{sec:sparse}. The proof relies on a fact that was recently used in several papers in
extremal graph theory, which states that graph of large enough average degree contain subgraphs with good expansion properties.
In the same section we also prove Theorem \ref{t:u5useful}.
Finally, in Section \ref{sec:two} we prove Theorem \ref{t:beta-2}.

\section{When $H$ is not a forest}\label{sec:prelim}

In this section we prove Propositions \ref{p:non-2-chrom}, \ref{p:2-chrom-upper} and \ref{p:non-strong}.
In what follows, when we refer to the transitive tournament $T_n$, we will assume, unless otherwise stated, that its vertices are $\{1,\ldots,n\}$
and $i < j$ implies that $(i,j)$ is an edge. We will frequently use the following well known fact (sometimes attributed to \cite{stearns-1959} and \cite{EM-1964})

\begin{theorem}\label{t:transitive}
Every tournament on $2^{h-1}$ vertices contains a copy of $T_h$.
\end{theorem}

\medskip

\Proof {\bf (of Proposition \ref{p:non-2-chrom}):} Let $H$ be a tournament with $h$ vertices and with $\chi(H)=r \ge 3$.
We first prove that $t(T_n,H) \ge \frac{r-2}{r-1}\binom{n}{2}$.
Consider the Tur\'an graph $T(n,r-1)$, that is the complete $(r-1)$-partite graph with $n$ vertices and with each color class
of size $\lfloor n/r \rfloor$ or $\lceil n/r \rceil$.
Add bidirectional edges to $T_n$ such that the undirected graph induced by the bidirectional
edges is isomorphic to $T(n,r-1)$. Observe that the resulting semi-complete digraph has chromatic number $r-1$ since each color class still induces a transitive set.
Hence, it does not contain $H$ as the latter has $\chi(H)=r$. We have added $|E(T(n,r-1))|$ edges to $T_n$ without introducing a copy of $H$.
Since $|E(T(n,r-1))| \ge \frac{r-2}{r-1}\binom{n}{2}$ we have that $t(T_n,H) \ge \frac{r-2}{r-1}\binom{n}{2}$.

We next prove that $t(n,H) \le \left(\frac{r-2}{r-1}+o_n(1)\right)\binom{n}{2}$.
Let the color classes of an $r$-coloring of $H$ be $H_1,\ldots,H_r$ where $|H_i|=h_i$ and notice that $h_i \le h-r+1$.
Let $U$ be the complete $r$-partite graph with $2^{h-r}$ vertices in each vertex class. By the Erd\H{o}s-Stone Theorem \cite{ES-1946},
$\mbox{ex}(n,U) = \left(\frac{r-2}{r-1}+o_n(1)\right)\binom{n}{2}$.
Given a tournament $G$, if we add to it $\mbox{ex}(n,U)$ edges, the subgraph on the bidirectional edges contains a copy of $U$.
Now consider the set $U_i$ of vertices of $G$ induced by the $i$'th color class of such a copy of $U$. It is a sub-tournament of $G$ on $2^{h-r} \ge 2^{h_i-1}$ vertices.
By Theorem \ref{t:transitive}, $U_i$ contains a subset $W_i$ which induces a copy of $T_{h_i}$.
Mapping $H_i$ (which induces a copy of $T_{h_i}$ in $H$) to $W_i$ and orienting all edges between $W_i$ and $W_j$ (recall that they are all bidirectional now)
as they are oriented between $H_i$ and $H_j$ in $H$, we obtain that $G$ together with its $\mbox{ex}(n,U)$ bidirectional edges contains a copy of $H$.
Hence, $t(n,H) \le \mbox{ex}(n,U)$ and the result follows.
\qed

\medskip

We now turn to the proof of Proposition \ref{p:2-chrom-upper}. Recall that for a $2$-chromatic tournament $H$, we used $s(H)$
to denote the smallest integer such that there is a $2$-coloring of $H$ with a
color class of size $s$. The other color class has $h-s \ge s$ vertices.
Denote the color classes of $H$ by $H_1$ and $H_2$ where $|H_1|=s$ and $|H_2|=h-s$.

\medskip

\Proof {\bf (of Proposition \ref{p:2-chrom-upper}):}
The proof is similar to the proof of the upper bound in Proposition \ref{p:non-2-chrom}.
Consider the complete bipartite undirected graph $U=K_{2^{s-1},2^{h-s-1}}$.
By the K\"ov\'ari-S\'os-Tur\'an Theorem \cite{KST-1954}, $\mbox{ex}(n,U)=O(n^{2-1/2^{s-1}})$.
If we take any tournament $G$ and add to it $\mbox{ex}(n,U)$ edges, the undirected graph formed by the bidirectional edges contains a copy of $U$.
Now consider the set $A$ of vertices of $G$ induced by the color class of such a copy of $U$ whose size is $2^{s-1}$ and the set $B$
induced by the color class of such a copy of $U$ whose size is $2^{h-s-1}$.
By Theorem \ref{t:transitive}, $A$ contains a subset $A'$ which induces a copy of $T_s$ and $B$ contains a subset $B'$ which induces a copy of $T_{h-s}$.
Mapping $H_1$ (which induces a copy of $T_s$ in $H$) to $A'$, and
mapping $H_2$ (which induces a copy of $T_{h-s}$ in $H$) to $B'$,
and orienting all edges between $A'$ and $B'$ (recall that they are all bidirectional now) as they are oriented between $H_1$ and $H_2$ in $H$,
we obtain that $G$ together with its $\mbox{ex}(n,U)$ bidirectional edges contains a copy of $H$.
Hence, $t(n,H) \le ex(n,U)$ and the result follows.
\qed

\medskip

For certain tournaments with $\chi(H) \leq 2$ one
can actually prove a much stronger bound than the one stated in Proposition \ref{p:2-chrom-upper}.
To this end we need to recall an interesting notion defined by Berger et al. \cite{BCCFLSST-2013}.
They say that a tournament $H$ a {\em hero} if there is a constant $c_H$ so that any tournament $T$ satisfying $\chi(T) > c_H$ contains a copy of $H$.
The main result of \cite{BCCFLSST-2013} is a precise characterization of heroes.
This characterization implies that all heroes have $\chi(H) \le 2$ and hence heroes form a particular interesting sub-family of $2$-chromatic tournaments.

\begin{theorem}\label{t:2-chrom-upper-hero}
If $H$ is a hero with $s=s(H)$, then $t(n,H)=O(n^{2-1/s})$\;.
\end{theorem}

\Proof
Set $U=K_{s,h-s}$. By the K\"ov\'ari-S\'os-Tur\'an Theorem \cite{KST-1954}, $\mbox{ex}(n,U)=O(n^{2-1/s})$.
Given a tournament $G$, if $G$ is not $H$-free, then $t(G,H)=0$\footnote{$t(G,H)$ is the least integer $t$ such that adding $t$ edges to $G$ results in a digraph containing a copy of $H$.}.
So assume that $G$ is $H$-free. Since $H$ is a hero, we know that $G$ is $c$-colorable, for some constant $c=c(H)$.
So let the color classes of G be $C_1,\ldots,C_c$. Now add to $G$ a set of $c^2 \cdot ex(n,U)$ bidirectional edges.
Now, if some $C_i$ contains $\mbox{ex}(n,U)$ bidirectional edges, then we are done.
This is because $C_i$ will then have a bidirectional copy of $U$. Since $C_i$ is transitive in $G$, the first class of such a copy induces in $G$ a transitive $T_s$
and the second class induces a transitive $T_{h-s}$. Since all edges between these classes are bidirectional, we obtain a copy of $H$ in $G[C_i]$ after adding the bidirectional edges.
If some pair $(C_i,C_j)$ contains $\mbox{ex}(n,U)$ bidirectional edges between $C_i$ and $C_j$, then we are done.
This is because $(C_i,C_j)$ will have a copy of $U$ where, without loss of generality, the first class of this copy is in $C_i$ and the second class is in $C_j$.
Since $C_i$ is transitive in $G$, the first class of such a copy induces in $G$ a transitive $T_s$
and the second class induces a transitive $T_{h-s}$. Since all edges between these classes are bidirectional, we obtain a copy of $H$ in $G[C_i \cup C_j]$ after adding the bidirectional edges.
Now, since there are only $\binom{c}{2}$ pairs $(C_i,C_j)$ and only $c$ classes $C_i$ and we have at least $(\binom{c}{2}+c)ex(n,U)$ bidirectional edges,
we must either have one of the two possibilities above occurring.
We have shown that $t(G,H) \le c^2 \cdot ex(n,U)$ so $t(n,H) \le O(n^{2-1/s})$, as claimed.
\qed

\medskip

We now turn to the proof of Proposition \ref{p:non-strong}.
We will actually prove the following more precise statement.

\begin{prop}\label{p:non-strong1}
Let $H$ be a tournament with $h$ vertices. If $H$ is not a forest, then:
$t(T_n,H)=\Omega(n^{1+\frac{4}{3h-\epsilon}})$ where
$\epsilon=4$ if $h \equiv 0 \bmod 4$,
$\epsilon=7$ if $h \equiv 1 \bmod 4$,
$\epsilon=6$ if $h \equiv 2 \bmod 4$,
$\epsilon=9$ if $h \equiv 3 \bmod 4$.
Moreover, for every $\epsilon > 0$, there exists $k$ and a tournament $H$ with a back-edge graph consisting of a matching of size $k$ which has $t(T_n,H)=\Omega(n^{2-\epsilon})$.
\end{prop}

\Proof We start with the first part of the proposition for which we need the following result.
\begin{lemma}\label{l:even-cycle}
For all $k \ge 2$ there are balanced bipartite graphs with $n$ vertices, with no cycles of length at most $2k$ and with
$\Omega(n^{1+\frac{2}{3k-3+\epsilon}})$ edges, where $\epsilon=0$ if $k$ is odd and $\epsilon=1$ if $k$ is even \cite{LUW-1995}.
\end{lemma}
Suppose $H$ has $h$ vertices, $\chi(H)=2$, but $H$ is {\em not} a forest tournament.
In other words, in any partition of $V(H)$ into two transitive sets $A$ and $B$, the set of edges going from $B$ to $A$ (viewed as an undirected bipartite graph) contains a cycle.
Now take $T_n$ and let $X$ be the first half of the vertices (vertices $1,\ldots, \lfloor n/2 \rfloor$) and $Y$ be the remaining vertices.
So all edges go from $X$ to $Y$. Add the maximum amount of bidirectional edges between $X$ and $Y$ while making sure that the set
of bidirectional edges (viewed as an undirected balanced bipartite graph with $n$ vertices) has no cycle of length $h$ or smaller.
Let $G$ denote the obtained semi-complete digraph. We claim that $G$ has no copy of $H$.
Indeed, suppose it had  a copy of $H$ and let $A \subset X$ be the vertices of such a copy falling in $X$ and $B \subset Y$ the vertices of the copy falling in $Y$.
Then this partitions $H$ into two transitive sets $A$ and $B$. But now, since $H$ is not a forest, the set of edges of this copy of $H$
going from $B$ to $A$ must contain a cycle, and in particular, a cycle of length at most $h$. But by construction, the bidirectional edges in $G$ are the only edges that
go from $B$ to $A$ and they do not contain a cycle of length $h$ or smaller, a contradiction.
Using Lemma \ref{l:even-cycle} with $k=\lfloor h/2 \rfloor$ we obtain that the number of edges added to $T_n$ to obtain $G$
is $\Omega(n^{1+\frac{4}{3h-\epsilon}})$ where $\epsilon=4$ if $h \equiv 0 \bmod 4$, $\epsilon=7$ if $h \equiv 1 \bmod 4$, $\epsilon=6$ if $h \equiv 2 \bmod 4$,
$\epsilon=9$ if $h \equiv 3 \bmod 4$. Hence, $t(T_n,H)$ is at least this large.

For the second part of the proposition, consider the tournament $\Delta_k$ having $k$ strongly connected components, each of which is the $3$-cycle $C_3$.
The components may be denoted by $Z_1,\ldots,Z_k$ where each $Z_i$ induces a $C_3$ on the vertices $(a_i,b_i,c_i)$ and all edges go from $Z_i$ to $Z_j$ when $i < j$.
An example of a back-edge graph of $\Delta_k$ consisting of just a matching of size $k$
is obtained by the order $a_1,b_1,c_1,a_2,b_2,c_2,\ldots,a_k,b_k,c_k$ where only the edges $(c_i,a_i)$ for $i=1,\ldots,k$ are back edges.
Observe that $\chi(\Delta_k)=2$ by taking one color class to be, say, $\{a_1,b_1,a_2,b_2,\ldots,a_k,b_k\}$ and the other class to be $\{c_1,\ldots,c_k\}$.
Nevertheless, it is not difficult to check that already $\Delta_3$ is not a tournament forest.
However, more can be said. Take any partition of the vertices of $\Delta_k$ into two transitive sets $A$ and $B$.
Then $Z_i$, being a $C_3$, must have at least one vertex in $A$ and at least one vertex in $B$.
In particular, every element of $B$ from $Z_i$ points t every element of $A$ from $Z_j$ where $j > i$.
Thus, there are at least $\binom{k}{2}$ edges pointing from $B$ to $A$.

Given $\epsilon > 0$, let $t=\lceil 2/\epsilon \rceil-1$. The well-known lower bound for complete bipartite Tur\'an numbers given in \cite{ES-1974} asserts
that for all $n$ sufficiently large, there are bipartite graphs (with $n/2$ vertices in each side) with $n^{2-2/(t+1)} \ge n^{2-\epsilon}$ edges that do not contain a copy of $K_{t,t}$.
Let $k$ be the least integer such that $\binom{k}{2} \ge 2(3k)^{2-1/t}$. By the K\"ov\'ari-S\'os-Tur\'an Theorem \cite{KST-1954}, any bipartite graph with $3k$ vertices
and more than $2(3k)^{2-1/t}$ edges contains a copy of $K_{t,t}$.
Now use the same construction and notation as in the proof of the first part of the proposition where now $H=\Delta_k$. We obtain a semi-complete digraph $G$ with at least
$n^{2-2/(t+1)} \ge n^{2-\epsilon}$ bidirectional edges where the bidirectional edges in $G$ are the only edges that
go from $B$ to $A$ and they do not contain a $K_{t,t}$, and hence there is no copy of $\Delta_k$. Thus, $t(T_n,\Delta_k) = \Omega(n^{2-\epsilon})$.
This proves the second part of the proposition.
\qed

\bigskip

\section{Tournaments vs matrices}\label{sec:equiv}

In this Section we prove Theorem \ref{t:ConjEquiv}. We start with the easy direction.

\begin{lemma}\label{l:equiveasy}
Conjecture \ref{c:PTC} implies Conjecture \ref{c:MainConj}.
\end{lemma}

\Proof
Suppose $H=(V,E)$ is a forest tournament and suppose it has a bipartition into two sets
$L=\{\ell_1,\ldots,\ell_{k_1}\}$ and $R=\{r_1,\ldots,r_{k_2}\}$ so that
the edges pointing from $R$ to $L$ form a forest (as in Definition \ref{TourForest}).
Let $M$ denote the $k_1 \times k_2$ matrix in which
$M_{i,j} =1$ if and only if $H$ has an edge pointing from $r_j$ to $\ell_i$. Since $H$ is a tournament
forest, then the matrix $M$ is a forest. Thus if Conjecture \ref{c:PTC} holds, then
$\mbox{ex}(n,M) = O(n(\log n)^p)$. We now show that if this is the case, then $t(T_n,H) =O( n(\log n)^{p+1})$.

Suppose we add edges to $T_n$ in such a way that for every integer $m$, the resulting digraph does not
contain two disjoint intervals of length $m$, with at least $m(\log m)^p$ bidirectional edges connecting them.
Then in particular the number of bidirectional edges connecting $\{1,\ldots,n/2\}$ to $\{n/2+1,\ldots,n\}$
is at most $n(\log n)^p$. Hence, if we let $f(m)$ denote the largest possible number of bidirectional edges within
an interval of length $m$, we get that $f(m) \leq 2f(m/2)+m(\log m)^p$ giving $f(m)=O(m(\log m)^{p+1})$.
Hence, if $G$ is a digraph resulting from adding $O(n(\log n)^{p+1})$ edges to $T_n$,
then we are guaranteed to have two intervals $X,Y$, with $X$ preceding $Y$, of length $m$ each,
with $m(\log m)^p$ bidirectional edges between them.
Suppose that $X=\{x_1,\ldots,x_m\}$ and $Y=\{y_1,\ldots,y_m\}$ and define
an $m \times m$ matrix $A$ by setting $A_{i,j}=1$ if and only if there is a bidirectional edge connecting
$y_j$ with $x_i$. Then $A$ contains at least $m(\log m)^{p+1}$ $1$'s and thus\footnote{Although we defined the problem
of bounding $\mbox{ex}(n,M)$ only with respect to square $k \times k$ matrices, it is clear that if Conjecture \ref{c:PTC} holds for
all forest $k \times k$ matrices then it holds also also for non-square matrices since we can just add rows/columns of $0$'s which
have negligible effect on $\mbox{ex}(n,M)$.} contains a copy of $M$.
It is now easy to see that if this copy of $M$ in $A$ has rows $i_1\ldots,i_k$ and columns $j_1,\ldots,j_k$, then
the vertices $x_{i_1},\ldots,x_{i_k},y_{j_1},\ldots,y_{j_k}$ span a copy of $H$ in $G$.
Hence if Conjecture \ref{c:PTC} holds, then so does Conjecture \ref{c:MainConj}.
\qed

\medskip

We now turn to prove that Conjecture \ref{c:MainConj} implies Conjecture \ref{c:PTC}, which will require some preparation.
Suppose $M$ is a $k \times k$ matrix with $0/1$ entries.
Note that in order to verify Conjecture \ref{c:PTC} it is enough to do so for matrices in which there is no all-zero row and no all-zero column.
Indeed, otherwise we can just add ones to such lines or columns while maintaining acyclicity.
Denote by $M(i,j)$ the entry in row $i$ and column $j$.

Given $M$ as above, let $M^*$ be the tournament defined on the vertex set $\{\ell_1,\ldots,\ell_k,r_1,\ldots,r_k\}$ and orient the edges as follows.
If $M(i,j)=1$, then we orient the edge from $r_j$ to $\ell_i$. If $M(i,j)=0$, then we orient the edge from $\ell_i$ to $r_j$.
We also orient from $\ell_i$ to $\ell_j$ and from $r_i$ to $r_j$ for all $1 \le i < j \le k$.
Notice that $M^*$ is a tournament forest since $M$ is acyclic. Its left side is $L=\{\ell_1,\ldots,\ell_k\}$ and its right side is $R=\{r_1,\ldots,r_k\}$
and the edges going from $R$ to $L$ form a forest.
Furthermore, every vertex of $L$ has a vertex from $R$ pointing to it and every vertex of $R$ points to a vertex of $L$.

Given $M^*$, as above, let $M^*_p$ be the tournament defined as follows.
Take $p$ copies of $L$ denoted $L_1\ldots,L_p$ where $L_s=\{\ell_{s,1},\ldots,\ell_{s,k}\}$,
and $p$ copies of $R$ denote $R_1,\ldots,R_p$, where $R_s=\{r_{s,1},\ldots,r_{s,k}\}$.
For each $1 \leq s \leq p$ put a copy of $M^*$ on $L_s \cup R_s$ by identifying\footnote{That is, the map sending
vertex $\ell_{s,i}$ of $M^*_k$ to vertex $\ell_i$ of $M^*$ and
vertex $r_{s,j}$ of $M^*_k$ to vertex $r_j$ of $M^*$
should be an isomorphism between $M^*_p[L_s \cup R_s]$ and $M^*$.} $\ell_{s,i}$ with $\ell_i$
and $r_{s,j}$ with $s_j$. Other than that, for every $s < t$ orient all edges between $L_s$ and $L_t$ to point from $L_s$
to $L_t$, all edges between $R_s$ and $R_t$ to point from $R_s$ to $R_t$, and all edges between $L_s$ and $R_t$ to point from $L_s$ to $R_t$.
Note that $M^*_p$ is also a tournament forest with partition $\cup_{s=1}^p L_s$ and $\cup_{s=1}^p R_s$.

Recall that $\mbox{ex}(n,M)$ is the least integer $m$ so that every $n \times n$ matrix $A$ with $m$ entries equal to $1$ contains the pattern matrix $M$.

\begin{lemma}\label{l:equivhard}
For every $M$ there is $p$ satisfying $$\mbox{ex}(n,M) \leq t(T_{2n},M_p^*)\;.$$
\end{lemma}

\Proof
Let $A$ be an $n \times n$ matrix with $\mbox{ex}(n,M)-1$ entries equal to $1$ and which does not contain the pattern $M$.
Let $T_{2n}$ be the transitive tournament on $\{1,\ldots,2n\}$, and
define a semi-complete digraph $G$ on these vertices by turning each pair $(i,j+n)$ with $1 \leq i \leq n$ and $1 \leq j \leq n$ into a bidirectional edge
if and only if $A(i,j)=1$. We will prove that $G$ has no copy of $M_p^*$ for all large enough $p$, and hence $t(T_{2n},M_p^*) \ge \mbox{ex}(n,M)$ as required.

Assume the contrary, and let a copy of $M^*_p$ in $G$ be spanned by a set of vertices $X$.
Recall that we denote the vertices of $M^*_p$
by $\ell_{s,i}$ and $r_{s,i}$ with $1\leq s \leq p$ and $1 \leq i \leq k$, where
$L_s=\{\ell_{s,1},\ldots,\ell_{s,k}\}$ and $R_s=\{r_{s,1},\ldots,r_{s,k}\}$.
Let $f:V(M^*_p) \mapsto X$ be an isomorphism from $M_p^*$ to the copy of $M^*_p$ in $G[X]$.
For each $1\leq s \leq p$ and $1 \leq i \leq l$ let $\ell'_{s,i}=f(\ell_{s,i})$ and $r'_{s,i}=f(r_{s,i})$.
For every $1 \leq s \leq p$ set $L'_s=\{\ell'_{s,1},\ldots,\ell'_{s,k}\}$ and $R'_s=\{r'_{s,1},\ldots,r'_{s,k}\}$.
Recalling the definition of $M^*_p$ this means that for every $1 \leq s \leq p$ the map sending $\ell'_{s,i}$ to $\ell_i$ and $r'_{s,i}$ to $r_i$
is an isomorphism between a subgraph of $G[L'_s \cup R'_s]$ and $M^*$.

Suppose first that for some $1 \leq s \leq p$ all the vertices of $L'_s$ appear (in $G$) before\footnote{Recall that the vertices of $G$ are
$\{1,\ldots,n\}$ so when we say ``before'' we refer to the natural order of the vertices.} those of $R'_s$.
Recall that by its construction, $G$ has no bidirectional edges within $\{1,\ldots,n\}$ and within $\{n+1,\ldots,2n\}$.
Hence, if $L'_s$ has a vertex in $\{n+1,\ldots,2n\}$, then this vertex has no edges pointing to it from $R'_s$.
This contradicts our assumption that $M$ has no zero rows.
For a similar reason $R'_s$ has no vertex in $\{1,\ldots,n\}$, as otherwise this contradicts our assumption that $M$ has no zero columns.
We see that $L'_s \subseteq \{1,\ldots,n\}$ and $R'_s \subseteq \{n+1,\ldots,2n\}$.
Since $L$ and $R$ span transitive sets in $M^*$ so should $L'_s$ and $R'_s$.
Since the only edges of $G$ within the sets $\{1,\ldots,n\}$ and $\{n+1,\ldots,2n\}$ are those of $T_{2n}$ we get that
the vertices of $L'_s$ and $R'_s$ appear in $G$ in the following order
$$
\ell'_{s,1} < \cdots < \ell'_{s,k} < r'_{s,1} < \cdots < r'_{s,k}\;.
$$
As noted above sending $\ell'_{s,i}$ to $\ell_i$ and $r'_{s,i}$ to $r_i$
is an isomorphism between a digraph spanned by $L'_s \cup R'_s$ and $M^*$.
This means that whenever $(r_i,\ell_j)$ is an edge of $M^*$ (or, equivalently, whenever $M(i,j)=1$) we must have a bidirectional edge between
$\ell'_{s,i}$ and $r'_{s,j}$.
But by the construction of $G$ from $A$ this means that we must have $A[\ell'_{s,i},r'_{s,j}-n]=1$.
Hence the sub-matrix of $A$ consisting of the rows $\{\ell'_{s,1},\ldots,\ell'_{s,k}\}$  and the columns $\{r'_{s,1}-n,\ldots,r'_{s,k}-n\}$
contains $M$ as a pattern, contradicting the assumption regarding $A$.

Suppose now that the condition in the previous paragraph does not hold.
Then for every $1 \leq s \leq p$ we have two vertices $r'_{s} \in R_s$ and $\ell'_{s} \in L_s$
so that $r'_{s}$ appears before $\ell'_{s}$.
Choose vertex $m$ in $G$ so that exactly $p/2$ of the above $p$ vertices $r'_{1},\ldots,r'_p$ appear before $m$.
Let $I \subseteq [p]$ be such that $s \in I$ if and only if $r'_{s}$ appears before $m$.
Then for every $s \not \in I$ we know that $\ell'_{s}$ appears after $m$.
We have thus found a set $X'$ of $p/2$ vertices $\{r'_{s}: s \in I\}$ and a set $Y'$ of $p/2$ vertices $\{\ell_{s}: s \not \in I\}$
so that all vertices of $X'$ appear in $G$ before all vertices of $Y'$. Since $M^*$ is a forest, the edges of
$M^*_k$ pointing from $\cup_{s=1}^p R_s$ to $\cup_{s=1}^p L_s$ form a forest. This means that for any pair of subsets $X \subseteq \cup_{s=1}^p R_s$
and $Y \subseteq \cup_{s=1}^p L_s$ of size $p/2$ each, there are less than $p$ edges pointing from $X$ to $Y$. Put differently,
there are at least $(p/2)^2-p+1$ edges pointing from $Y$ to $X$. Since we assume
that $f$ is an isomorphism between the a subgraph of $G[X]$ and $M^*_p$ and $X' \subseteq \cup_{s=1}^p R'_s$ and $Y \subseteq \cup_{s=1}^p L'_s$
this means that in $G$ we must have at least $(p/2)^2-p+1$ edges pointing from $Y'$ to $X'$. But since, $X'$ appears before $Y'$
this means that we have at least $(p/2)^2-p+1$ bidirectional edges between $X'$ and $Y'$.
This implies that in the sub-matrix of $A$ consisting of the $p/2$ rows $Y'$ and the $p/2$ columns\footnote{Recall that the vertices of
$G$ are the integers $\{1,\ldots,2n\}$ so by $X'-n$ we mean the set of columns whose indices are $\{x'-n~:~x' \in X'\}$.} $X'-n$
we have  at least $(p/2)^2-p+1$ entries equal to $1$. By the K\"ovari-Sos-Tur\'an Theorem, every $(p/2) \times (p/2)$ matrix with $k(p/2)^{2-1/k}$ entries
equal to $1$, contains a $k \times k$ sub-matrix all of whose entries are $1$. Hence, for any $p$ satisfying $(p/2)^2-p+1 > k(p/2)^{2-1/k}$, the
matrix $A$ will contain a copy of $M$ as a pattern, which again contradicts the assumption.
\qed

\medskip

\Proof {\bf (of Theorem \ref{t:ConjEquiv}):} Follows from Lemmas \ref{l:equiveasy} and \ref{l:equivhard}
and the discussion preceding Lemma \ref{l:equivhard}.
\qed

\medskip

\Proof {\bf (of Theorem \ref{t:non-linear}):}
Consider the matrices in Figure \ref{f:1}.
F\"uredi and Hajnal \cite{FH-1992} showed that the matrix $M_1$ has $\mbox{ex}(n,M_1)=\Theta(n \log n)$.
Lemma \ref{l:equivhard} thus proves that $t(T_{n},(M_1^*)_p) = \Omega(n \log n)$ for $p$ sufficiently large (in fact, it suffices to take $p=40$ in this case since
a bipartite graph with $20$ vertices in each side and at least $400-39=361$ edges contains a $K_{3,3}$).
Pettie \cite{pettie-2010} showed that the matrix $M_2$ has $\mbox{ex}(n,M_2)=\Omega(n \log n \log \log n)$.
Lemma \ref{l:equivhard} thus proves that $t(T_{n},(M_2^*)_p) = \Omega(n \log n \log\log n)$ for $p$ sufficiently large (in fact, he proved that already the rectangular matrix $M_2$
with the last row deleted satisfied this lower bound).
\qed

\begin{figure}[ht]
$$
\begin{array}{cc}
M_1 = \left[\begin{array}{ccc} \bullet &  \bullet & \\  \bullet &  &  \\ & & \bullet \end{array}\right] &
M_2 = \left[\begin{array}{ccccc} & \bullet & & \bullet & \bullet \\  &  &\bullet & &\bullet \\& \bullet & & & \\\bullet&  & & & \bullet\\ & & & & \bullet\end{array}\right]
\end{array}
$$
\caption{The matrix $M_1$ with $\mbox{ex}(n,M_1)=\Theta(n \log n)$ and the matrix $M_2$ with $\mbox{ex}(n,M_2)=\Omega(n \log n \log\log n)$.
Bullets represent $1$ and blanks represent $0$.}\label{f:1}
\end{figure}

\section{Bounding $t(n,H)$ using sparse families}\label{sec:sparse}

The main result we prove in this section is Theorem \ref{t:sparse} which reduced
the problem of bounding $t(n,H)$ to finding ``parse'' characterizations'' of prime $H$-free tournaments.
We will end this section with the proof of Theorem \ref{t:u5useful} showing that the tournament $U_5$
has such a sparse family. We will need the following two lemmas, whose
statements use the notions introduced in Subsection \ref{subsec:PTCapp}.

\begin{lemma}\label{l:findexp}
If $U$ is an undirected $n$-vertex graph with $bn(\log n)^{b+1}$ edges,
then $U$ has an $m$-vertex subgraph which is a $b\log^bm$-expander.
\end{lemma}

\Proof  Consider the following process in which we iteratively construct a sequence of graphs
$U_0,U_1,\ldots$, with the property that if $U_i$ has $m_i$ vertices than it has at least
$bm_i(\log m_i)^{b+1}$ edges. We start by setting $U_0=U$. Assuming we already defined
$U_i$ we do the following. If $U_i$ is a $b\log^bm_i$-expander we stop. Otherwise
we do the following; let $X$ be a set of size $x \leq m_i/2$ having less than $bx\log^{b}m_i$
edges connecting it to its complement (within $U_i$). If the graph spanned by $X$ has at least
$bx\log^{b+1}x$ edges, then we set $U_{i+1}$ to be this graph. Otherwise,
if the graph spanned by the complement of $X$ has at least
$b(m_i-x)(\log(m_i-x))^{b+1}$ edges, then we set $U_{i+1}$ to be this graph.
We claim that it cannot be the case that both cases fail.
Indeed, if this is the case, then the total number of edges in $U_i$ is at most
$$
bx(\log (m_i/2))^{b+1}+b(m_i-x)\log^{b+1}m_i+bx\log^{b}m_i < bm_i\log^{b+1} m_i
$$
contradicting the assumption that $U_i$ has at least $bm_i(\log m_i)^{b+1}$ edges.
Hence the above process must eventually stop
with the last graph $U_i$ being an $m_i$-vertex $b\log^bm_i$-expander.
\qed

\medskip

\begin{lemma}\label{l:expnonprime}
If $G=(V,E)$ is a semi-complete $n$-vertex digraph which is a $4\log n$-expander, then with probability at least $3/4$
a tournament drawn from ${\cal D}_G$ is not prime.
\end{lemma}

\Proof Suppose $T \sim{\cal D}_G$. Take any set $X \subseteq V$ with $x=|X|$ where $2 \le x \le n-1$.
We will prove that the probability that $X$ is homogenous in $T$ is at most
$(2n \binom{n}{x})^{-1}$. We consider first the case $|X| \le n/2$.
Let $K$ be the set of bidirectional edges of $G$ connecting $X$ and $V \setminus X$.
Since $G$ is assumed to be a $4\log n$-expander we know that $|K| \ge 4 x \log n$.
Consider some $w \in V \setminus X$. We say that $w$ is Type-$1$ if there is at least one non-bidirectional edge connecting $w$ and $X$ (or, in other words, not all $x$ edges between $w$ and $X$ are bidirectional edges). Otherwise, $w$ is Type-$2$.
Let $K_w \subset K$ be the bidirectional edges incident with $w$. Now, if $w$ is Type-$1$, there is some non-bidirectional edge $e$ connecting $w$ and $X$,
so in order for $X$ to be homogenous in $T$, it must be the case the for each bidirectional edge $\{w,i\} \in  K_w$ the orientation that was chosen
for $\{w,i\}$ must agree with the (fixed) orientation of $e$.
This occurs with probability $(1/2)^{|K_w|}$.
Suppose now that $w$ is Type-$2$. In this case, have that $|K_w|=x \ge 2$.
In a similar manner, we get that conditioning on any orientation of one of the edges of $K_w$,
the probability that all the other edges get the same orientation is $(1/2)^{|K_w|-1} \leq (1/2)^{|K_w|/2}$.
We see that in both cases, the probability that $w$ would have the property that all pairs ${w,j}$ with $j \in X$
with have the same orientation is at most $(1/2)^{|K_w|/2}$.
So overall, the probability that $X$ is homogenous is at most
$$
\prod_{w \in V \setminus X} 2^{-|K_w|/2}= 2^{-\frac12\sum_{w \in V \setminus X}|K_w|} = 2^{-\frac12|K|} \le 2^{-2x \log n} < \frac{1}{4n\binom{n}{x}}\;.
$$
Notice that if $|X| \ge n/2$, then the same argument holds replacing $x$ with $n-x$ and using the fact that $\binom{n}{x}=\binom{n}{n-x}$.
We have thus proved that the probability that $T$ contains a homogenous set is at most
$$
\sum_{X \subset V\;,\;2 \le |X| \le n-1} \frac{1}{4n\binom{n}{x}} \le \frac{1}{4}\;.
$$
Hence with probability at least $3/4$, the tournament $T$ is prime. \qed

\medskip

The proof of Theorem \ref{t:sparse} now follows quite easily from the above lemmas.

\medskip

\Proof {\bf (of Theorem \ref{t:sparse})}: Suppose there is a sparse $H$-useful family of tournaments ${\cal Q}$
and suppose $b$ is the constant from Definition \ref{def:sparse}.
Set $t=\max\{b,4\}$. We claim that $t(n,H) \leq tn(\log n)^{t+1}$. To see this,
take any $n$-vertex tournament and add $tn(\log n)^{t+1}$ edges to it.
Let $G=(V,E)$ be the resulting semi-complete digraph. Suppose to the contrary that $G$ is $H$-free.
Let $U$ be the undirected graph on the same vertex set $V$ in which $(i,j)$ is an edge if and only if
$\{i,j\}$ is a bidirectional edge in $G$. Then $U$ contains $tn(\log n)^{t+1}$ edges, so we can use Lemma \ref{l:findexp}
in order to find an $m$-vertex subgraph $U'$ of $U$ which is a $t\log^{t}m$-expander.
Suppose $U'$ has vertex set $V'$
and let $G'=(V',E')$ be the subgraph of $G$ induced by $V'$.
Let us now pick a tournament $T \sim {\cal D}_{G'}$. Since $t \geq 4$, Lemma \ref{l:expnonprime} tells us that the probability
that $T$ is prime is at most $1/4$ and the assumption that ${\cal Q}$ is sparse together with the fact that $t \geq b$
implies that the probability that $T \in {\cal Q}$ is at most $1/2$. Hence, with positive probability, $T$ will be
a prime and will not belong to ${\cal Q}$. But since $G$ is assumed to be $H$-free, then so must be $G'$ and $T$.
We have thus found an $H$-free prime tournament $T$ which does not belong to ${\cal Q}$, contradicting the assumption
that $Q$ is $H$-useful.  \qed

\medskip

We end this section with the proof that the family ${\cal Q}_5$ (defined in Subsection \ref{subsec:PTCapp}) is sparse.

\medskip

\Proof {\bf (of Theorem \ref{t:u5useful}):} Suppose $G$ is a $10\log n$-expander. Observe that for odd $n$, the tournament $C_n$
is $(n-1)/2$-regular. For a given vertex $v$, let $d_{b}(v)$, $d_{in}(v)$, $d_{out}(v)$ be the number of bidirectional edges
touching $v$, the number of edges pointing at $v$ and the number of edges pointing from $v$.
Suppose that from $k$ of the $d_{b}(v)$ bidirectional edges touching $v$ the edge
pointing from $v$ was chosen, and from the other $d_{b}(v)-k$ bidirectional edges, the edge pointing to $v$ was chosen.
Then for $v$ to have out-degree $(n-1)/2$ we must have $d_{out}(v)+k=d_{in}(v)+d_{b}(v)-k$.
Thus, setting $d=d_{b}(v)$ and $d'=(d_{in}(v)+d_{b}(v)-d_{out}(v))/2$, the probability of this event is at most $\binom{d}{d'}/2^d$.
Since $G$ is a $\log n$-expander, we have\footnote{Note that here we are only using the fact that a in $\log n$-expander, every
 vertex has degree at least $\log n$.} $d \geq \log n$, implying that
this probability is $O(1/\sqrt{\log n})$, implying that the expected
number of vertices of degree $(n-1)/2$ is $O(n/\sqrt{\log n})$ so the probability that they will all have
degree $(n-1)/2$ is $o(1)$. In particular, the probability that $T \sim {\cal D}_G$ will be isomorphic to $C_n$
is $o(1)$.

We now show that the probability that $T \sim {\cal D}_G$ will satisfy the second condition of Theorem \ref{t:prime}
is less than $1/3$. Consider a given partition of $n$ into three nonnegative integers $x_1,x_2,x_3$ such that $x_1+x_2+x_3=n$, $x_1 \le x_2 \le x_3$
and a partition of $V(G)$ into three sets $V_1, V_2, V_3$ of sizes $|V_i|=x_i$.
We estimate the probability that $T \sim {\cal D}_G$ will have the property that $T[V_1 \cup V_2]$, $T[V_1 \cup V_2]$, and $T[V_2 \cup V_3]$ are transitive. Denote this property by ${\cal P}$.

Assume first that $x_2=0$. In this case $x_3=n$ and $T$ will have ${\cal P}$ if and only if $T$ is transitive. Since $G$ is a $10 \log n$-expander, it has
at least $5n\log n$ bidirectional edges, hence the probability that the direction of these edges will comply with any of the $n!$ possible
ordering of $T$'s vertices is at most $n!2^{-5n \log n} \ll 1/n^2$.

Assume therefore that $x_2 \ge 1$. Then there are at least $10x_2\log n$ bidirectional edges between $V_2$ and $V_1 \cup V_3$, so there are at least $5x_2\log n$ bidirectional edges
between $V_2$ and $V_j$ for one of $j=1,3$. The probability that $V_2 \cup V_j$ is transitive in $T$ is at most
$$
\binom{x_2+x_j}{x_2}2^{-5x_2\log n} < \binom{n}{x_2}2^{-5x_2\log n} < 2^{-4x_2 \log n} < \frac{1}{n^{x_1+x_2+2}} \;.
$$
The number of ways to choose sets $V_1,V_2,V_3$ with the given sizes $x_1,x_2,x_3$ is $\binom{n}{x_1}\binom{n-x_1}{x_2} < n^{x_1+x_2}$.
Hence, for a given partition $x_1+x_2+x_3=n$, the probability that some partition of $V(G)$ into three sets will cause $T$ to satisfy ${\cal P}$
is less than $1/n^2$. As there are less than $n^2/4$ ways to make a partition of $n$ to $x_1+x_2+x_3$, we get by the union bound that the probability
that $T \sim {\cal D}_G$ satisfies the second condition of Theorem \ref{t:prime} is at most $1/4+1/n^2 < 1/3$.

Altogether, the probability that $T \sim {\cal D}_G \in {\cal Q}_5$ is less than $1/2$, so ${\cal Q}_5$ is sparse.
\qed

\section{Tournaments with at most two back edges}\label{sec:two}

We prove Theorem \ref{t:beta-2}.
Suppose that $H$ is a tournament with a back-edge graph consisting of two edges.
We may assume that both edges have no common endpoint as otherwise $H$ is a star tournament and we are done.
Suppose therefore that the back edges are $(a,b),(c,d)$ where $a,b,c,d$ are distinct. There are three possible configurations that need to be handled.
The {\em disjoint} configuration, where both $a,b$ precede $c,d$ (namely $b < a < d < c$), the {\em intersecting} configuration where $b < d < a < c$ and the
{\em containment} configuration, where $b < d < c < a$.
Lemmas \ref{l:disjoint}, \ref{l:containment} and \ref{l:intersecting} handle each of these configurations.
Before that we need to prove the following useful lemma.
Suppose $G$ is an undirected ordered graph with vertices $\{1,\ldots,n\}$.
A subgraph $G^*$ of $G$ is a {\em $t$-gap} if any pair of vertices $u,v$ of $G^*$ satisfy $|u-v| \ge t$. Furthermore, vertices $1,\ldots,t-1$ and vertices $n-t+2,\ldots,n$
are not in $G^*$.
\begin{lemma}\label{l:gap}
Suppose that $G$ is an undirected ordered graph with vertices $\{1,\ldots,n\}$ and with $m$ edges.
Then $G$ has a $t$-gap subgraph with at least $(m-3nt)/(et)^2$ edges.
\end{lemma}
\Proof
For $i=t,\ldots,n-t+1$, vertex $i$ is selected according to the following rule.
If any one of the vertices $i-1,i-2,i-t+1$ was selected, $i$ is never selected. Otherwise, $i$ is selected with probability $1/t$.
This guarantees that the set of selected vertices induce a $t$-gap subgraph.
Now, suppose $u$ and $v$ are vertices with $v-u \ge t$ where $u \ge t$ and $v \le n-t+1$.
Observe that $u$ is selected with probability at least $(1-1/t)^{t-1}\cdot(1/t) \ge 1/(et)$. Now, $v$ is selected with probability at least $1/(et)$ regardless of whether
$u$ is selected or not. Hence, both $u$ and $v$ are selected with probability at least $1/(et)^2$ (this probability can be slightly improved with a somewhat more
detailed analysis). There are less than $3nt$ edges $(u,v)$ that are incident with $1,\ldots,t$ or with $n-t+2,\ldots,n$ or have $|u-v| < t$.
So, if $G^*$ is the subgraph of $G$ induced by the set of chosen vertices, we have that the expected number of edges of $G^*$ is
at least $(m-3nt)/(et)^2$, as required.
\qed

\begin{lemma}\label{l:disjoint}
Suppose that $H$ is a tournament on $h$ vertices with an order having a back-edge graph consisting of two edges in a disjoint configuration.
Then $t(T_n,H) \le 20h^2n$.
\end{lemma}
\Proof
Let $X_h$ denote the undirected tree on the set of vertices $\{x_1,\ldots,x_h\}$ and whose edges are
$(x_1,x_{h-1}),(x_2,x_{h-1}),\ldots,(x_{h-2},x_{h-1}),(x_{h-2},x_h)$.

Consider $T_n$ and the semi-complete digraph obtained after adding to it $2hn$ edges.
The set of bidirectional edges is an ordered undirected graph $G$ on $\{1,\ldots,n\}$.
We claim that $G$ either has two edges $(x_j,x_i),(x_\ell,x_k)$ with $x_i < x_j < x_k < x_\ell$ or else it contains an ordered copy of $X_h$.

The claim is proved by induction on $n$ where the case $h < n$ holds trivially.
We classify the edges of $G$ as follows. If both endpoints are in the first half of vertices (vertices $1,\ldots,\lfloor n/2 \rfloor$), the edge is called {\em left}.
If both are in the second half of vertices, it is called {\em right}. Otherwise, it is called {\em crossing}.
Now, if there is some left edge and also some right edge, then we are done as the left edge can play the role of $(x_j,x_i)$ and the right edge can play the role of $(x_\ell,x_k)$.
Otherwise, if there are $2h\lfloor n/2 \rfloor$ left edges we are done by induction on $T_{\lfloor n/2 \rfloor}$.
Otherwise, if there are $2h\lceil n/2 \rceil$ right edges we are done by induction on $T_{\lceil n/2 \rceil}$.
So we are left with the case that there are at least $2h\lfloor n/2 \rfloor \ge (h-1)n$ crossing edges.
Hence, there is a subgraph $G^*$ of $G$ consisting only of crossing edges and with minimum degree at least $h-1$.
Consider the shortest edge of $G^*$, namely an edge $(u,v)$ where $u \in \{1,\ldots,\lfloor n/2 \rfloor\}$ and $v \in \{\lfloor n/2 \rfloor+1,\ldots,n\}$ and $v-u$ is as small as possible.
Let $x_1,\ldots,x_{h-2}$ be $h-2$ distinct neighbors of $v$ in $G^*$ with $x_1 < x_2 < \cdots < x_{h-2} < u$.
Let $x_h$ be some neighbor of $u$ other than $v$ and notice that $x_h > v$ since $v-u$ is minimal. Setting $u=x_{h-2}$ and $v=x_{h-1}$ we obtain that $G^*$, and therefore also $G$,
contains an ordered copy of $X_h$.

Having proved this claim we now proceed to the proof of the statement of the lemma. So consider $T_n$ and
the semi-complete digraph obtained after adding to it $20h^2n$ edges. Let $G$ be the ordered graph on the vertices of $T_n$ and the $20h^2 n$ bidirectional edges.
By Lemma \ref{l:gap}, $G$ has an $h$-gap subgraph $G'$ with a least $(20h^2n - 3hn)/(eh)^2 \ge 2n$ edges.
Notice also that $G'$ has at most $n/h$ vertices. By the claim in the previous paragraph, $G'$ either has two
edges $(x_j,x_i),(x_\ell,x_k)$ with $x_i < x_j < x_k < x_\ell$ or else it contains an ordered copy of $X_h$.
Consider now the tournament $H$ with vertex set $\{y_1,\ldots,y_h\}$ and with the order $y_1,\ldots,y_h$ having two back edges.
Suppose these two back edges are $(y_q,y_p)$ and $(y_s,y_r)$ where by the assumption on them being in disjoint configuration we have $p < q < r < s$.

Consider first the case that $G'$ has two edges $(x_j,x_i),(x_\ell,x_k)$ with $x_i < x_j < x_k < x_\ell$.
Then we can find a copy of $H$ in $T_n \cup G$ by
assigning $y_1,\ldots,y_{p-1}$ to $p-1$ vertices smaller than $x_i$,
assigning $y_p$ to $x_i$,
assigning $y_{p+1},\ldots,y_{q-1}$ to $q-p-1$ vertices larger than $x_i$ and smaller than $x_j$,
assigning $y_q$ to $x_j$,
assigning $y_{q+1},\ldots,y_{r-1}$ to $r-q-1$ vertices larger than $x_j$ and smaller than $x_k$,
assigning $y_r$ to $x_k$,
assigning $y_{r+1},\ldots,y_{s-1}$ to $s-r-1$ vertices larger than $x_k$ and smaller than $x_\ell$,
assigning $y_s$ to $x_\ell$,
and assigning $y_{s+1},\ldots,y_h$ to $h-s$ vertices larger than $x_\ell$.

We may now assume that $G'$ has an ordered copy of $X_h$.
Consider the order of $H$ obtained by placing $y_p$ after $y_r$, hence the order
$y_1,\ldots,y_{p-1},y_{p+1},\ldots,y_r,y_p,y_{r+1},\ldots,y_h$.
The back-edge graph of this order is isomorphic to $X_{r-p+1}$, as $y_p$ now has back edges to $y_{p+1},\ldots,y_{q-1},y_{q+1},\ldots,y_r$ and $y_r$ gets a back edge from $y_s$.
Since $G'$ has a copy of $X_h$ and since $r-p+1 \le h$ we have in particular that $G'$ has a copy of $X_{r-p+1}$.
Since $G'$ is an $h$-gap, we can find a copy of $H$ in $T_n \cup G$ realizing this order of $H$ as follows:
assigning $y_1,\ldots,y_{p-1}$ to $p-1$ vertices smaller than $x_1$ (recall that the vertices of $X_{r-p-1}$ are $x_1,\ldots,x_{r-p+1}$),
assigning $y_{p+1},\ldots,y_{q-1}$ to vertices $x_1,\ldots, x_{q-p-1}$,
assigning $y_q$ to a vertex larger than $x_{q-p-1}$ and smaller than $x_{q-p}$,
assigning $y_{q+1},\ldots,y_r$ to vertices $x_{q-p},\ldots,x_{r-p-1}$,
assigning $y_p$ to $x_{r-p}$,
assigning $y_{r+1},\ldots,y_{s-1}$ to vertices larger than $x_{r-p}$ and smaller than $x_{r-p+1}$,
assigning $y_s$ to $x_{r-p+1}$,
and assigning $y_{s+1},\ldots,y_h$ to $h-s$ vertices larger than $x_{r-p+1}$.
\qed

\begin{lemma}\label{l:containment}
Suppose that $H$ is a tournament on $h$ vertices with an order having a back-edge graph consisting of two edges in a containment configuration.
Then $t(T_n,H) \le 30hn$.
\end{lemma}
\Proof
Consider $T_n$ and the semi-complete digraph obtained after adding to it $3n+1$ edges.
Let $G$ be the ordered graph on the vertices of $T_n$ and the $3n+1$ bidirectional edges.
We claim that $G$ has two edges $(x_\ell,x_i),(x_k,x_j)$ with $x_i < x_j < x_k < x_\ell$.

Let the {\em forward} degree of a vertex $u$ of $G$, denoted by $f(u)$, be the number of edges $(u,v)$ with $v > u$.
Notice that $\sum_{i=1}^n f(i)=3n+1$. Let $U$ be the set of vertices with $f(u) > 0$.
For $u \in U$ let $r(u)$ be the rightmost neighbor of $u$ and let $pred(u)$ be the predecessor of $u$ in $U$ (if $u$ is the first vertex of $U$,
then $pred(u)=\phi$ and also define $r(\phi)=0$).
Let $u$ be the least vertex in $U$ such that $r(u) \le r(pred(u))+f(u)-2$.
Notice that this initially does not occur as for the smallest vertex in $u \in U$ we have $r(u) \ge u+f(u)$.
But also notice that $u$ must exists, as otherwise we would have for the maximal element $v \in U$ that
$r(v) \ge \sum_{u \in U} (f(u)-2)$ but this cannot happen as $r(v) \le n$ while $\sum_{u \in U} (f(u)-2) \ge (3n+1)-2n > n$.
So, we indeed have some vertex $u$ such that $r(u) \le r(pred(u))+f(u)-2$.
But if the rightmost neighbor of $u$ is in location $r(u)$, then there is also a forward neighbor of $u$, say $w$, with $w \le r(u)-(f(u)-1)$.
But then we have $pred(u) < u < w < r(pred(u))$ so letting  $x_i=pred(u)$, $x_j=u$, $x_k=w$ and $x_\ell=r(pred(u))$ the claim follows.

Having proved this claim we now proceed to the proof of the statement of the lemma. So consider $T_n$ and
the semi-complete digraph obtained after adding to it $30hn$ edges.
Let $G$ be the ordered graph on the vertices of $T_n$ and the $30h n$ bidirectional edges.
By Lemma \ref{l:gap}, $G$ has an $h$-gap subgraph $G'$ with a least $(30hn - 3hn)/(eh)^2 \ge 3n/h+1$ edges.
Notice also that $G'$ has at most $n/h$ vertices.
By the claim in the previous paragraph, $G'$ has two edges $(x_\ell,x_i),(x_k,x_j)$ with $x_i < x_j < x_k < x_\ell$.
Consider now the tournament $H$ with vertex set $\{y_1,\ldots,y_h\}$ and with the order $y_1,\ldots,y_h$ having two back edges.
Suppose these two back edges are $(y_s,y_p)$ and $(y_r,y_q)$ where by the assumption on them being in containment configuration we have $p < q < r < s$.
We can find a copy of $H$ in $T_n \cup G$ by
assigning $y_1,\ldots,y_{p-1}$ to $p-1$ vertices smaller than $x_i$,
assigning $y_p$ to $x_i$,
assigning $y_{p+1},\ldots,y_{q-1}$ to $q-p-1$ vertices larger than $x_i$ and smaller than $x_j$,
assigning $y_q$ to $x_j$,
assigning $y_{q+1},\ldots,y_{r-1}$ to $r-q-1$ vertices larger than $x_j$ and smaller than $x_k$,
assigning $y_r$ to $x_k$,
assigning $y_{r+1},\ldots,y_{s-1}$ to $s-r-1$ vertices larger than $x_k$ and smaller than $x_\ell$,
assigning $y_s$ to $x_\ell$,
and assigning $y_{s+1},\ldots,y_h$ to $h-s$ vertices larger than $x_\ell$.
\qed

\begin{lemma}\label{l:intersecting}
Suppose that $H$ is a tournament on $h$ vertices with an order having a back-edge graph consisting of two edges in an intersecting configuration.
Then $t(T_n,H) \le 20hn$.
\end{lemma}
\Proof
Consider $T_n$ and the semi-complete digraph obtained after adding to it $2n-2$ edges.
Let $G$ be the ordered graph on the vertices of $T_n$ and the $2n-2$ bidirectional edges.
We claim that $G$ has two edges $(x_k,x_i),(x_\ell,x_j)$ with $x_i < x_j < x_k < x_\ell$.
This immediately from the fact that the maximum number of edges in an outer-planar graph with $n$ vertices is $2n-3$.

Having proved this claim we now proceed to the proof of the statement of the lemma. So consider $T_n$ and
the semi-complete digraph obtained after adding to it $20hn$ edges.
Let $G$ be the ordered graph on the vertices of $T_n$ and the $20h n$ bidirectional edges.
By Lemma \ref{l:gap}, $G$ has an $h$-gap subgraph $G'$ with a least $(20hn - 3hn)/(eh)^2 \ge 2n/h-2$ edges.
Notice also that $G'$ has at most $n/h$ vertices.
By the claim in the previous paragraph, $G'$ has two edges $(x_k,x_i),(x_\ell,x_j)$ with $x_i < x_j < x_k < x_\ell$.
Consider now the tournament $H$ with vertex set $\{y_1,\ldots,y_h\}$ and with the order $y_1,\ldots,y_h$ having two back edges.
Suppose these two back edges are $(y_r,y_p)$ and $(y_s,y_q)$ where by the assumption on them being in intersecting configuration we have $p < q < r < s$.
We can find a copy of $H$ in $T_n \cup G$ by using the same embedding as in Lemma \ref{l:containment}. \qed

\medskip

The proof of Theorem \ref{t:beta-2} now follows immediately from
Lemmas \ref{l:disjoint}, \ref{l:containment} and \ref{l:intersecting}. \qed

\bibliographystyle{plain}

\end{document}